\documentclass[12pt,reqno]{amsart}
\usepackage{verbatim}
\usepackage{amsthm}
\usepackage{amsmath}
\usepackage{amssymb}
\usepackage{mathrsfs}
\usepackage{pdfsync, graphicx}
\usepackage{color}

\begin{document}

\numberwithin{equation}{section}
\renewcommand{\theequation}{\thesection.\arabic{equation}}
\newcommand{\cblue}{\color{blue}}
\newcommand{\cred}{\color{red}}
\newcommand{\cviolet}{\color{violet}}

\newcommand{\Ac}{\mathcal A}
\newcommand{\Asc}{\mathscr A}
\newcommand{\Bc}{\mathcal B}
\newcommand{\Bsc}{\mathscr B}
\newcommand{\Ec}{\mathcal E}
\newcommand{\Ectau}{\mathcal{E}_{\tau}}
\newcommand{\Ect}{\tilde{\Ec}}
\newcommand{\Esc}{\mathscr E}
\newcommand{\eps}{\epsilon}
\newcommand{\Fc}{\mathcal F}
\newcommand{\Fsc}{\mathscr{F}}
\newcommand{\Fsco}{\mathscr{F}_0}
\newcommand{\Gamk}{\Gamma_k}
\newcommand{\Gsc}{\mathscr{G}}
\newcommand{\Gc}{\mathcal{G}}
\newcommand{\Gsm}{\mathcal{G}_{\lambda, \mu}}
\newcommand{\Ic}{\mathcal{I}}
\newcommand{\Isc}{\mathscr{I}}
\newcommand{\Lsc}{\mathscr{L}}
\newcommand{\Mcl}{{\mathcal M} \, }
\newcommand{\Msc}{\mathscr M}
\newcommand{\Mcltau}{\mathscr{M}_{\tau}}
\newcommand{\N}{\mathbb{N}}
\newcommand{\Nb}{\mbox{$\mathbb N$}}
\newcommand{\Om}{\Omega}
\newcommand{\bdy}{\partial\Omega}
\newcommand{\Omb}{\overline{\Omega}}
\newcommand{\pal}{\partial}
\newcommand{\Pc}{\mathcal P}
\newcommand{\Qc}{\mathcal Q} 
\newcommand{\R}{\mathbb R}
\newcommand{\Rlam}{\R_{\lambda_k}^2}
\newcommand{\Rc}{\mathcal R}
\newcommand{\Rb}{\overline{\R}}
\newcommand{\Rn}{{\R}^n}
\newcommand{\RN}{{\R}^N}
\newcommand{\Rp}{[0,\infty)}
\newcommand{\ra}{\rightarrow}
\newcommand{\Sc}{\mathcal S}
\newcommand{\Sctau}{\Sc_{\tau}}
\newcommand{\ub}{\bar{u}}
\newcommand{\uh}{\hat{u}}
\newcommand{\ut}{\tilde{u}}
\newcommand{\vt}{\tilde{v}}


\newcommand{\aip}[1]{\langle#1\rangle}					
\newcommand{\mtau}{m_{\tau}}
\newcommand{\n}[1]{\left\vert#1\right\vert}
\newcommand{\nm}[1]{\left\Vert#1\right\Vert}
\newcommand{\ipb}[2]{[#1, #2]_b}
\newcommand{\ipmo}[2]{[#1, #2]_{m_1}}						
\newcommand{\ipmt}[2]{[#1, #2]_{m_2}}
\newcommand{\ipm}[2]{[#1, #2]_m}
\newcommand{\nmb}[1]{\left\Vert #1\right\Vert_b}
\newcommand{\nmmo}[1]{\left\Vert #1\right\Vert_{m_1}}
\newcommand{\nmmt}[1]{\left\Vert #1\right\Vert_{m_2}}
\newcommand{\nmm}[1]{\left\Vert #1 \right\Vert_m}


\newcommand{\Div}{\mathop{\rm div}\nolimits}
\newcommand{\dsg}{d \sigma}
\newcommand{\gradu}{\nabla u}
\newcommand{\gradv }{\nabla v}
\newcommand{\Gscmu}{\Gsc(.,\mu)}
\newcommand{\Iby}{ \int_{\bdy}} 
\newcommand{\IOm}{\int_{\Om}}

\newcommand{\s}{\varsigma}
\newcommand{\lam}{\lambda}
\newcommand{\lamo}{\lambda_{(1)}}
\newcommand{\lamt}{\lambda_{(2)}}
\newcommand{\muo}{\mu_{(1)}}
\newcommand{\mut}{\mu_{(2)}}

\newcommand{\Lamtau}{\Lambda_{\tau}}
\newcommand{\wuo}{\hat{w}_{1}}
\newcommand{\wut}{\hat{w}_{2}}

\newcommand{\Elamu}{E_{\lambda, \mu}}
\newcommand{\Jlamu}{J_{\lambda, \mu}}
\newcommand{\Ilamu}{I_{\lambda, \mu}}
\newcommand{\rlamu}{r_{\lambda, \mu}}


\newcommand{\barr}{\begin{eqnarray}}
\newcommand{\beq}{\begin{equation}}
\newcommand{\bpf}{\begin{proof} \quad}
\newcommand{\btm}{\begin{thm}}
\newcommand{\blem}{\begin{lem}}
\newcommand{\bcor}{\begin{cor}}
\newcommand{\ecor}{\end{cor}}
\newcommand{\elem}{\end{lem}}
\newcommand{\earr}{\end{eqnarray}}
\newcommand{\eeq}{\end{equation}}
\newcommand{\epf}{\end{proof}}
\newcommand{\etm}{\end{thm}}


\newcommand{\Hone}{H^1(\Om)}
\newcommand{\Hones}{\Hone \times \Hone} 
\newcommand{\Hz}{H_{0}^1(\Om)}
\newcommand{\Linb}{L^{\infty}(\bdy,\dsg)}
\newcommand{\Lp}{L^p (\Om)}
\newcommand{\Lpb}{L^p (\bdy, \dsg)}
\newcommand{\Lq}{L^q (\Om)}
\newcommand{\Lt}{L^2 (\Om)}
\newcommand{\Ltb}{L^2 (\bdy,\dsg)}
\newcommand{\RxHone}{\R \times \Hone}
\newcommand{\Vp}{V_{+}}
\newcommand{\Vz}{V_{0}}
\newcommand{\Vk}{V_{k}}
\newcommand{\Wk}{W_{k}}
\newcommand{\Wktd}{\widetilde{W}_k}

%
\newtheorem{thm}{Theorem}[section]
\newtheorem{cor}[thm]{Corollary}
\newtheorem{cond}{Condition}
\newtheorem{lem}[thm]{Lemma}
\newtheorem{prop}[thm]{Proposition}
%
%


\title[Eigencurves and Elliptic Operators]{Eigencurves for linear elliptic equations}

\author[Rivas \& Robinson]{M. A. Rivas \, and \, Stephen B. Robinson}

\address{Mauricio A. Rivas \newline
Department of Mathematics and Statistics, Wake Forest University \newline
PO Box 7388, 127 Manchester Hall,
Winston-Salem, NC 27109, USA}
\email{rivasma@wfu.edu} 

\address{Stephen B. Robinson \newline
Department of Mathematics and Statistics, Wake Forest University \newline
PO Box 7388, 127 Manchester Hall,
Winston-Salem, NC 27109, USA}
\email{sbr@wfu.edu}

\date{\today}

\begin{abstract} 
This paper provides results for eigencurves associated with self-adjoint linear elliptic boundary value problems.
The elliptic problems are treated as a general two-parameter eigenproblem for a triple $(a, b, m)$ of continuous symmetric bilinear forms on a real separable Hilbert space $V$.
Variational characterizations of the eigencurves associated with $(a, b, m)$ are given and various orthogonality results for corresponding eigenspaces are found.
Continuity and differentiability properties, as well as asymptotic results, for these eigencurves are proved.
These results are then used to provide a geometrical description of the eigencurves.
\end{abstract}
\maketitle


\section{Introduction}


This paper is motivated by the study of eigencurves associated with self-adjoint linear elliptic boundary value problems such as
\beq\label{e.robin1}
\begin{aligned}
-\Delta u(x) \; = & \; \mu\, m_0(x)\, u(x) \quad\text{for }x\in \Omega \\
D_{\nu} u(x) \; + \; c(x)\,u(x)\; = & \; \lambda\, b_0(x)\, u(x)\quad \text{for }x\in \partial\Omega,
\end{aligned}
\eeq
where $c, b_0, m_0$ are given functions in appropriate $L^p$-spaces on a bounded region $\Omega$ of $\R^N$ satisfying mild boundary regularity requirements, and $\lambda, \mu$ are real eigenparameters.
Here, $m_0$ is assumed to be strictly positive, $b_0$ may be sign-changing, and $D_{\nu}$ denotes the outer normal derivative.
For such problems, the boundary and interior equations may be combined in weak form using bilinear forms.
Therefore, our focus in this article is on the analysis of abstract eigencurve problems associated with triples $(a, b, m)$ of continuous symmetric bilinear forms on a real Hilbert space $V$.

Our main result generalizes the geometric characterization of eigencurves given for {\it Sturm-Liouville} problems in Binding and Volkmer \cite{BV96}.
We expand on the issues treated in \cite{BV96} regarding continuity, differentiability, and asymptotics of eigencurves, and also provide results for issues not appearing in the ODE case.

The analysis in this paper is based on the use of spectral results for bilinear forms obtained in Auchmuty \cite{Au4}, and in some respects our work may be regarded as a complementary continuation of that paper.
The use of bilinear forms provides a simpler alternative to the usual operator-theoretic approach to eigencurves as it avoids the use of dual, or other, Sobolev spaces and operators.
For a treatment of eigenproblems using bilinear forms invoking properties of associated linear operators see Attouch, Buttazzo, and Michaille \cite{ABM}, or Blanchard and  Br\"{u}ning \cite{BB92}.  
For a classical treatment of eigencurves using the theory of closed operators on Hilbert space see Kato \cite{K}.

A paper on the use of eigencurves to establish existence results for some indefinite weight semilinear elliptic problems is Ko and Brown \cite{KB}; those problems arise, for instance, in population genetics.
Their results are based on results for the principle eigencurve(s) for linear boundary value problems with indefinite weight and Robin boundary conditions given in Afrouzi and Brown \cite{AB99}.

The papers of Mavinga and Nkashama \cite{MN2010} and Mavinga \cite{Mavinga} establish existence results for some nonlinear elliptic equations with nonlinear boundary equations where the nonlinearities interact (in some sense) with the associated generalized Steklov-Robin spectrum.  The spectra considered in \cite{MN2010}, \cite{Mavinga} may be regarded as sections (slices) of eigencurves associated to their equations.  
It is worth noting that 
the principle as well as {\it higher} eigenvalues of the related linear problems are considered in those papers.  

Recently, two-parameter problems for the Laplacian have been used in Section 10 of Auchmuty and Rivas \cite{AR16} to provide representations of Sobolev spaces on product regions as {\it tensor products} of Sobolev spaces on the individual factor regions.
That paper forms a contemporary reference for the analysis of problems on product regions arising, for instance, in fluid mechanics, electromagnetic theory and elsewhere.

For a nice summary on other applications of eigencurves and connections to indefinite inner product spaces, see Binding and Volkmer \cite{BV96}.

The current paper is organized as follows.
The weak form of problem \eqref{e.robin1} and other concrete examples are discussed in the next section.  
These are merely applications used for expository purposes and to indicate the diversity of problems that may be treated using our general bilinear form framework.

The specific two-parameter eigenvalue problem for triples $(a, b, m)$ of bilinear forms considered in this paper is described in Section \ref{s.definitions}, as well as associated terminology and definitions to be used.
Section \ref{s.vcurves} details the variational characterizations for the eigencurves that are based on the constructive algorithm given in \cite{Au4}.
An immediate consequence of these characterizations is the concavity of the first eigencurve.

Special inner products are used in Section \ref{s.ortho} to establish orthogonality of certain eigenspaces associated with the eigencurves.
Then in Section \ref{s.continuity}, continuity of each eigencurve is established. 
In particular, the eigencurves are shown to be Lipschitz continuous.

We show in Section \ref{s.differentiability} that eigencuves are differentiable except possibly where they intersect.  
At an intersecting point, an eigencurve may have a corner, but has well-defined one-sided derivatives.  
Explicit formulae for one-sided derivatives at each point $\lambda\in \mathbb{R}$ are found in terms of spectral data for the pair $(b, m)$ of bilinear forms.

Another difficulty encountered here is in describing the asymptotic behaviour of eigencurves as the $b, m$ forms are not necessarily variants of each other.
For the Robin-Steklov problem \eqref{e.RS} (the weak form of problem \eqref{e.robin1}), the bilinear form $b$ comprises a boundary integral whereas the $m$ form comprises an interior-region integral as in \eqref{e.bmforms}, and thus a comparison of these forms is not straightforward.
However, it is shown in Section \ref{s.asymptotic} that the asymptotic behaviour of eigencurves is governed by the spectrum of $(b, m)$ and that the sign of the quadratic form $u\mapsto b(u, u)$ plays a role in these matters.

An interesting side issue encountered in our analysis is the possible appearance of straight lines within the {\it spectrum} (the collection of graphs of eigencurves) of $(a, b, m)$ due to the degeneracy of $b$.  
Degeneracy may occur in elliptic problems when weight functions are zero on a set of positive measure or when the equations involve boundary integrals;  for the Robin-Steklov problem \eqref{e.RS}, the function $b_0(x)$ may be zero on a set of positive $(N-1)$-dimensional Hausdorff measure, but the form itself is already degenerate as it comprises a boundary integral defined on $\Hone$.
Results for straight lines within the spectrum are given in Section \ref{s.independence} and are shown to be related to a question of linear independence of certain collections of eigenvectors.
Some simple examples of triples $(a, b, m)$ are provided in this section to illustrate the results.

In Section \ref{s.geometry} we prove a generalization of a main result in \cite{BV96} stating that any straight line intersects the first $n$ eigencurves 
in at most $2n$ points.
The straight line is assumed not a subset of the spectrum of $(a, b, m)$.
It is worth noting that even in the most degenerate case, there are at most countably many horizontal lines in the spectrum.
For our main result, we focus on the number of connected components the graph of an eigencurve may have above the intersecting line instead of considering points of intersections as eigenvalues in our general setting may have multiplicity greater than one.
This approach combined with the results of previous sections yields our generalized geometrical characterization of eigencurves.

\section{Motivating examples}\label{s.examples}

The abstract results of this paper apply to very general linear elliptic two-parameter eigenproblems.
In this section, we describe  a few examples that motivated our analysis, and we will use them to simplify the presentation and illustrate the general results of this paper.  
The examples include the (weak form of the) Robin-Steklov problem, the Sturm-Liouville problem treated in \cite{BV96}, and a generalized matrix eigenvalue problem.

In discussing the Robin-Steklov problem, the terminology and notation of Evans and Gariepy \cite{EG}  will be used here except that $\sigma, d\sigma$ will denote Hausdorff $(N-1)$-dimensional measure and integration with respect to this measure, respectively.
The real Lebesgue spaces $\Lp$ and $\Lpb$, $1\leq p\leq \infty$, are defined in the standard manner with the usual $p$-norms denoted $\|u\|_{p, \Om}$ and $\|u\|_{p,\partial \Om}$, respectively.
The gradient of a function $u$ will be denoted $\nabla u$, and $\Hone$ is the usual real Sobolev space with the standard $H^1$-inner product
\[
[u, v]_{1, 2} : = \IOm [ \nabla u \cdot \nabla v \, + \, u\, v]\, dx \quad \text{and corresponding norm}\quad \|u\|_{1,2} \; : = \; \IOm [|\nabla u|^2 + u^2|\, dx.
\]
When treating the Robin-Steklov problem \eqref{e.RS} below it is required that the embedding of $\Hone$ into $\Lp$ (for $1\leq p< p_s$ where $p_s = 2N/(N-2)$ if $N\geq 3$ and $p_s = \infty$ if $N=2$) and the trace operator $\Gamma:\Hone\rightarrow \Ltb$ are both compact.
This requirement holds, for instance, when $\Omega\subset \R^N$ is a bounded domain such that its boundary consists of a finite number of closed Lipschitz surfaces of finite surface area (see Chapter 4 of \cite{EG} for details on these issues).
This requirement is assumed throughout this paper.

The {\it two-parameter Robin-Steklov eigenproblem} is the problem of finding $(\lambda, \mu)\in \R^2$ such that there is a nonzero $u\in \Hone$ satisfying
\beq\label{e.RS}
\int_{\Omega}\nabla u\cdot \nabla v\, dx \, + \, \int_{\partial\Omega}c\, u\, v\, d\sigma \; = \; 
\lambda\, \int_{\partial \Omega} b_0\, u\, v \, d\sigma \, + \, \mu \, \int_{\Omega} m_0\, u\, v\, dx \qquad \text{for all }v\in \Hone.
\eeq
We shall assume $b_0, c\in \Linb$, $m_0\in \Lp$ with $p>N/2$, is given data satisfying $c\geq 0$ on $\partial\Omega$ and $\|c\|_{1,\partial\Omega} \, >0$, and $m_0$ is positive on $\Om$ with $\|m_0\|_{1, \Omega}\, >0$.  Here, $\lambda, \mu$ are real eigenparameters, and we emphasize that $b_0$ may be sign-changing.
The bilinear forms on $\Hone$ associated with \eqref{e.RS} are then 
\beq
a(u,v) \; : =  \; \int_{\Omega}\nabla u\cdot \nabla v \, dx \, + \, \int_{\partial\Omega}c\, u\, v \, d\sigma\, ,
\eeq
and
\beq\label{e.bmforms}
b(u, v)\; : =  \int_{\partial\Omega} b_0\, u\, v\, d\sigma \, \qquad \qquad \text{and} \qquad \qquad m(u, v) \; : = \: \int_{\Omega}m_0\, u\, v\, dx.
\eeq
The {\it two-parameter Robin-Steklov eigenpoint equation} \eqref{e.RS} then becomes:
\beq
a(u, v) \; = \; \lambda\, b(u, v) \, + \, \mu\, m(u, v) \qquad \text{for all }v\in \Hone.
\eeq 
It follows from Section 7 and 8 of \cite{Au4} that $a, b, m$ satisfy (A1)-(A3) given in Section \ref{s.definitions} below, so that the results of this paper apply to this Robin-Steklov eigenproblem.

The {\it Sturm-Liouville} eigenproblem given in \cite{BV96} is that of finding nontrivial $y=y(t)$ satisfying 
\beq \label{e.SL}
-(p(t)y'\, )'  +  q(t)\, y \; = \; (\lambda\, r(t)  +  \mu)\, y\qquad\text{on an interval }[t_0, t_1] \subset \mathbb{R},
\eeq
and the separated boundary conditions
\beq\label{e.SLBC}
\cos(\alpha)y(t_0) \, - \, \sin(\alpha)p(t_0)\, y'(t_0) \; = \; 0, \quad \qquad \cos(\beta)y(t_1) - \sin(\beta)p(t_1)\, y'(t_1) \; = \; 0.
\eeq
It is assumed in \cite{BV96} that $p$ is continuously differentiable and positive on $[t_0, t_1]$, $q$ and $r$ are piecewise continuous on $[t_0, t_1]$ and $\alpha, \beta \in \mathbb{R}$;
the eigenparameters are also $\lambda, \mu$.
It is is shown in \cite{BV96} that there is a sequence $\{\mu_n(\lambda) :n\in \mathbb{N}\}$ of simple analytic curves with variational characterization, that $\mu_1(\lambda)$ is concave with positive eigenfunctions, and that the curves satisfy a line intersection property.  That is, if $\ell$ is  a fixed straight line in $\mathbb{R}^2$, then $\ell$ intersects the first $n$ curves at most $2n$ times.


Lastly, the {\it two-parameter generalized matrix eigenproblem} is the problem of finding $(\lambda,\mu)\in \mathbb{R}^2$ and a nonzero vector $x\in \mathbb{R}^N$ satisfying
\beq \label{e.matrixproblem}
Ax \; = \; \lambda\, Bx \, + \, \mu\, x 
\eeq
where $A$ is a symmetric positive definite matrix, such as the Laplacian matrix, and $B$ a symmetric matrix on $\mathbb{R}^N$.
It is a good non-trivial exercise to show that this problem has a sequence of continuous curves $\{\mu_n(\lambda) : n\in \mathbb{N}\}$ with variational characterization, with the first eigencurve $\mu_1(\lambda)$ concave but not necessarily simple and isolated.

\begin{figure}[h]
\begin{center}
\includegraphics[width=110mm, height=75mm]{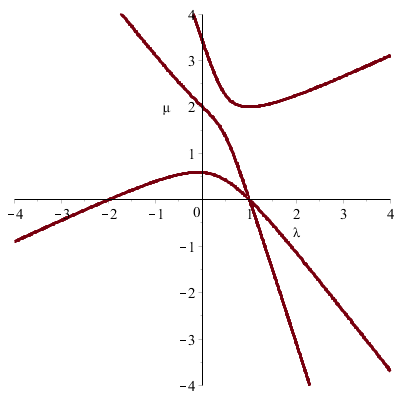}
\caption{Eigencurves for the matrix eigenproblem \eqref{e.matrixproblem}-\eqref{e.matrixvalues}.}
\label{matrix}
\end{center}
\end{figure}

When $A, B$ are taken to be the $3\times 3$ matrices given by
\beq\label{e.matrixvalues}
A\; = \; \begin{bmatrix}
2 & -1 & 0 \\
-1 & 2 & -1 \\
0 & -1 & 2
\end{bmatrix}
\qquad \text{and}\qquad 
B\; = \; \begin{bmatrix}
2 & -1 & 0 \\
-1 & 2 & -1 \\
0 & -1 & 0
\end{bmatrix}
\eeq
the eigencurves may be found (implicitly) and are plotted in Figure \ref{matrix}.

\section{Two-parameter eigenvalue problems for bilinear forms}\label{s.definitions}

The two-parameter eigenproblem will be studied in the framework of bilinear forms on Hilbert space with the following definitions and notation.
$V$ will denote a real, separable, Hilbert space with inner product and norm on $V$ denoted by $\langle \cdot, \cdot\rangle$ and $\|\cdot\|$, respectively.

Our interest is in describing the pairs $(\lambda, \mu)\in \R^2$ for which there is a nonzero $u$ in $V$ satisfying
\beq\label{e.eigenequation}
a(u, v) \; = \; \lambda\, b(u, v) \; + \; \mu\, m(u, v)  \qquad \text{for all }v\in V,
\eeq
where $a, b, m$ are bilinear forms on $V$ subject to the conditions described below.  
This will be called the $(a, b, m)$-{\it eigenproblem}. 
The pair $(\lambda, \mu)$ is said to be an {\it eigenpoint} of $(a, b, m)$ if there is a nonzero vector $u\in V$ satisfying \eqref{e.eigenequation}, and such $u$ will be called an {\it eigenvector} of $(a, b, m)$ corresponding to $(\lambda, \mu)$.  
The subset $\sigma(a, b, m)$ of $\R^2$ consisting of all eigenpoints $(\lambda, \mu)$ will be called the {\it spectrum} of $(a, b, m)$, and \eqref{e.eigenequation} will be called the {\it eigenpoint equation}.

Define $\Ac, \Bc, \Mcl$ to be the quadratic forms on $V$ associated with $a, b, m$ so that
\beq\label{e.qform}
\Ac(u) \; := \; a(u, u)\, , \qquad \qquad \Bc(u) \; : = \; b(u, u)\, , \qquad \text{and}\qquad \Mcl(u) \; := \; m(u, u).
\eeq
The assumptions on the bilinear forms to be used in this paper will include

\noindent
{\bf (A1):}\quad $a(\cdot, \cdot)$ is a continuous, symmetric, bilinear form on $V$ that is also $V$-coercive.
That is, there are constants $0< \kappa_1 \leq  \kappa_2 < \infty$ such that
\beq
\kappa_1 \, \|u\|^2 \; \leq \; \Ac(u) \; \leq \; \kappa_2\, \|u\|^2 \qquad \text{for all }u\in V.
\eeq
{\bf (A2):}\quad $b(\cdot, \cdot)$ is a weakly continuous, symmetric bilinear form on $V$.\\
{\bf (A3):}\quad  $m(\cdot, \cdot)$ is a weakly continuous, symmetric, bilinear form that satisfies
\[
\Mcl(u)\; > \; 0 \quad\text{for all nonzero }u\in V.
\]

When (A1) holds, then $a(\cdot, \cdot)$ defines an inner product on $V$ equivalent to the $V$-inner product and is called the $a$-{\it inner product}; the associated norm will be denoted by $\|\cdot\|_a$.
When (A3) holds, the quadratic form $\Mcl$ is strictly positive on $V$ so that $m$ is an inner product on $V$ and the  associated norm will be denoted by $\|u\|_m$.  
A vector $u$ in $V$ will be said to be $m$-{\it normalized} provided $m(u, u) =1$ .  
A bilinear form $b$ satisfing $(A2)$ may be negative or zero for some $v\in V$, and is said to be an {\it indefinite form} whenever $\Bc$ attains both positive and negative values.
 
When working merely with a pair of bilinear forms, the standard $(a, m)$-{\it eigenproblem} for a pair $(a, m)$ of symmetric, continuous bilinear forms on $V$ is that of finding $\mu \in \R$ and nontrivial $u\in V$ satisfying
\beq\label{e.standardeigen}
a(u, v) \; = \; \mu \, m(u, v) \qquad \text{for all } v\in V.
\eeq
A nonzero $u$ satisfying \eqref{e.standardeigen} is called an {\it eigenvector} associated to the {\it eigenvalue} $\mu$ of $(a, m)$.  
The number of linearly independent eigenvectors of $(a, m)$ corresponding to $\mu$ is called the {\it multiplicity} of $\mu$.    
When the multiplicity of $\mu$ is one, then $\mu$ is said to be a {\it simple} eigenvalue.  
The set $\sigma(a, m)$ of all distinct eigenvalues of $(a, m)$ will be called the {\it spectrum} of $(a, m)$.

When $(A1)$, $(A3)$ hold, there is an increasing sequence $\Lambda_0 :=\{\mu_{0n} : n\in \N\}$ of strictly positive eigenvalues of $(a, m)$, repeated according to (finite) multiplicity with $\mu_{0n} \rightarrow \infty$ as $n\rightarrow \infty$, 
and there is an associated sequence $\Ec_0 :=\{e_{0n} :n\in \N\}$ of eigenvectors constituting a basis for $V$ that is orthogonal with respect to both $a$ and $m$.  
These results follow from the analysis of Section 4 of \cite{Au4} and the notation here is chosen to simplify the presentation. 


\section{Variational characterization of eigencurves}\label{s.vcurves}

In this section the variational eigencurves associated to $(a, b, m)$ to be studied in this paper are obtained using the constructive algorithm given in Auchmuty \cite{Au4}.  
The construction yields a one-parameter family of sequences of eigenpoints for $(a, b, m)$ as well as a corresponding family of sequences of eigenvectors.

For $\lambda, \tau \in \R$ define $a_{\lambda, \tau}$ to be the bilinear form on $V$ given by
\beq\label{e.alamtau}
a_{\lambda, \tau}(u, v) \; : = \; a(u, v) \; - \; \lambda \, b(u,v) \; + \; \tau\, m(u, v) \qquad \text{for }u, v\in V.
\eeq
The basic coercivity result for these forms is the following.

\btm\label{t.tau}
Assume $a, b, m$ satisfy $(A1)$-$(A3)$ and $\Ac$ is defined by \eqref{e.qform}.
Then for fixed $r_0>0$ there exists a constant $\tau> 0$ such that for each $\lambda\in [-r_0, \, r_0]$ the bilinear form $a_{\lambda, \tau}$ given by \eqref{e.alamtau} satisfies 
\beq
a_{\lambda, \tau}(u, u)  \; \geq \; \frac{1}{2} \Ac(u) \qquad \text{for all }u\in V.
\eeq
\etm

\bpf
Arguing by contradiction, suppose there are sequences $\tau_n>0$, $\lambda_n\in [-r_0, \, r_0]$ and $u_n\in V$, with $\tau_n$ increasing to $+\infty$, satisfying
\beq\label{e.coercive1}
a_{\lambda_n, \tau_n}(u_n, u_n)\;  < \; \frac{1}{2}\, \Ac(u_n) \qquad \text{for all }n\in \N.
\eeq
The form $a_{\lambda_n, \tau_n}(u_n, u_n)$ is equal to $\Ac(u_n) - \lambda_n\Bc(u_n) + \tau_n\Mcl(u_n)$ where  $\Ac, \Bc, \Mcl$ are given by \eqref{e.qform}.
Without loss of generality, suppose $\lambda_n \rightarrow \lambda$ in $[-r_0, \, r_0]$, the $u_n$ are $a$-normalized so that $\Ac(u_n) =1$.
By assumption (A1) the set $\{u\in V: \Ac(u)=1\}$ is closed and bounded so we may also assume without loss of generality that $u_n$ converges weakly to $u$ in $V$.
Then \eqref{e.coercive1} becomes
\beq\label{e.coercive2}
\frac{1}{2} \; - \; \lambda_n\, \Bc(u_n) \; + \; \tau_n \Mcl(u_n) \; < \; 0 \qquad \text{for all }n\in \N.
\eeq
By weak continuity of $\Mcl, \Bc$ we see that $\Mcl(u_n) \rightarrow \Mcl(u)$ and $\Bc(u_n) \rightarrow \Bc(u)$.
If $u\neq0$, then $\tau_n \, \Mcl(u_n)\rightarrow +\infty$ which contradicts \eqref{e.coercive2}.  
If, on the other hand, $u=0$, then $\lambda_n \Bc(u_n) \rightarrow 0$ which implies
\[
\frac{1}{2} \; - \; \lambda_n \, \Bc(u_n) \; + \; \tau_n\, \Mcl(u_n) \; > \; 0
\]
for large enough $n$.  This again contradicts \eqref{e.coercive2}, so the assertion of the theorem  holds.
\epf

This results says that for $\lambda\in \R$ fixed, $a_{\lambda, \tau}$ will be a coercive bilinear form whenever $\tau$ is large enough and provided (A1), (A2) also hold. 
For the remainder of this section suppose $\lambda\in \R$ and $\tau>0$ have been chosen to satisfy this criteria.
Then $a_{\lambda, \tau}$ defines an inner product on $V$ equivalent to the $V$-inner product.

Denote by $\mu_1(\lambda)+\tau$ the first eigenvalue of $(a_{\lambda,\tau}\, ,  m)$.
It may be found by maximizing the weakly continuous functional $\Mcl$ defined by \eqref{e.qform} on the closed convex subset $C_{1, \lambda, \tau}$ of $V$ defined by $C_{1,\lambda, \tau} :=\{u\in V : a_{\lambda, \tau}(u,u) \leq 1\}$.  
The maximizers of this problem are eigenvectors of $(a_{\lambda, \tau}\,,  m)$ corresponding to the eigenvalue $\mu_1(\lambda)+\tau$
whereas the value of this problem is 
\beq
\frac{1}{\mu_1(\lambda) +\tau} \; = \; \sup_{u\in C_{1, \lambda, \tau}} \Mcl(u).
\eeq
These results follow directly from Section 3 of \cite{Au4} with the appropriate notational modifications.
Then using the construction of Section 4 of that paper an infinite sequence of eigenvalues and eigenvectors of $(a_{\lambda, \tau}, m)$ may be found.

Let $\mu_n(\lambda) +\tau >0$, with $n=1, 2, \ldots, k$, be the first $k$ successive smallest eigenvalues of $(a_{\lambda, \tau}\, , m)$ and let $e_{\lambda n}$ be associated $a_{\lambda, \tau}$-orthonormal eigenvectors.  
The $(k+1)^{st}$ eigenvector of $(a_{\lambda, \tau}\, , m)$ will be a maximizer of $\Mcl$ over the subset
\beq
C_{k+1, \lambda, \tau} \; := \; \{u\in C_{1, \lambda, \tau} \; :\; a_{\lambda, \tau}(u, e_{\lambda n})=0 \; \text{ for }n=1, 2, \ldots, k\}
\eeq
and the value of this $(k+1)^{st}$ problem is
\beq
\frac{1}{\mu_{k+1}(\lambda)+ \tau} \; = \; \sup_{u\in C_{k+1, \lambda, \tau}}\Mcl(u).
\eeq
Let $\Lambda_{\lambda} : =\{\mu_n(\lambda) \, : \, n\in \N\}$ be the sequence of such values repeated according to multiplicity and in increasing order and 
$\Ec_{\lambda} : =\{e_{\lambda n} \, :\, n\in\N\}$ be a sequence of associated eigenvectors constructed by this iterative process. 
The existence and some properties of this eigendata is described in Section 4 of \cite{Au4}.  
Specifically, the analysis there yields the following result.

\btm\label{t.abm.basis}
Assume $a, b, m$ satisfy $(A1)$-$(A3)$, $a_{\lambda, \tau}$ is defined by \eqref{e.alamtau} for $\tau$ large enough, and $\Lambda_{\lambda}, \, \Ec_{\lambda}$ are the sequences defined as above.  Then each $\mu_n(\lambda)$ in $\Lambda_{\lambda}$ is of finite multiplicity with $-\infty < \mu_1(\lambda) \leq \mu_2(\lambda) \leq \cdots \rightarrow +\infty$ and $\Ec_{\lambda}$ is an $a_{\lambda, \tau}$-orthonormal basis of $V$ consisting of vectors satisfying
\beq\label{e.eigeneq}
a(e_{\lambda n}\, , v) \; =\; \lambda\, b(e_{\lambda n}\, , v)  \; + \; \mu_n(\lambda) \, m(e_{\lambda n}\, , v) \qquad \text{for all }v\in V.
\eeq
\etm

\bpf
Conditions $(A1), (A4)$ of \cite{Au4} hold with the bilinear forms $a, m$ there taken to be the pair $(a_{\lambda, \tau}\, , m)$ in this paper, so the analysis of Section 4 there is applicable.  
This yields the sequences $\Lambda_{\lambda}$ and $\Ec_{\lambda}$ for the pair $(a_{\lambda, \tau}, m)$ and their properties.  
Rearranging the corresponding eigenequation for $(a_{\lambda, \tau}, m)$ gives the eigenpoint equation \eqref{e.eigeneq},
so that the $e_{\lambda n}$ are indeed eigenvectors of $(a, b, m)$ corresponding to eigenpoints $(\lambda, \mu_n(\lambda))$.
\epf

Define the $n^{th}$-{\it variational eigencurve} associated to $(a, b, m)$ to be the graph of the function $\lambda \mapsto \mu_n(\lambda)$.
This will be a subset of $\sigma(a, b, m)$ in $\R^2$, and we shall often say $\mu_n(\lambda)$ is the $n^{th}$ eigencurve associated to $(a, b, m)$.
The point $(\lambda_*, \mu_*) \in \sigma(a, b, m)$ is said to be an {\it eigenpoint of multiplicity $k_0$ for} $(a, b, m)$ accordingly as the value $\mu_*$ is an eigenvalue of multiplicity $k_0$ for the pair $(a_{\lambda_*b}\, , m)$ where $a_{\lambda_*b} : = a - \lambda_*b$.

When $a, b, m$ satisfy (A1)-(A3) and $\lambda\in \R$ is fixed, the following {\it minimax} characterization of eigenvalues will be useful for the analysis of these variational eigencurves:
\beq\label{e.minmax}
\mu_k(\lambda) \; \; = \; \; \inf_{Z_k}\sup_{\substack{v\in Z_k \\ v\neq 0}} \left\{ \frac{\Ac(v) \, - \, \lambda\, \Bc(v)}{\Mcl(v)}  \right\}
\eeq
where the infimum is taken over all subspaces $Z_k$ of $V$ of dimension $k$.
The supremum in \eqref{e.minmax} is equal to $\mu_k(\lambda)$ when $Z_k$ is a subspace generated by the first $k$ $m$-orthonormal eigenvectors corresponding to the smallest eigenvalues of $(a_{\lambda b}\, , m)$ where $a_{\lambda b}\ : = a - \lambda b$.  
These results follow from Section 4 of \cite{Au4} with $(a_{\lambda, \tau}\, , m)$ in place of $(a, m)$ and rearranging the relations obtained there to eliminate the parameter $\tau$.

A direct consequence of this characterization is the following.

\blem\label{l.mu1concave}
Assume $a, b, m$ satisfy $(A1)$-$(A3)$.  Then the first eigencurve $\mu_1(\lambda)$ is concave.
\elem

\bpf
It follows from \eqref{e.minmax} that $\mu_1(\lambda)$ is the infimum of affine functions on $\R$.
\epf


\section{Orthogonality relations and eigenspaces}\label{s.ortho}

Here some {\it orthogonality relations} among eigenspaces of $(a, b, m)$ will be described.  
In particular, distinct eigenspaces are shown to be orthogonal with respect to various bilinear forms whenever the corresponding eigenpoints lie on the same horizontal line or the same vertical line.
These results simplify the analysis of eigencurves.

Two vectors $v, w\in V$ are said to be orthogonal with respect to a bilinear form $\tilde{a}$ provided $\tilde{a}(v, w)=0$.
The forms $\tilde{a}$ arising here are, in general, indefinite on the whole space $V$ and are therefore not necessarily inner products on $V$.
A first result is the following.

\btm\label{t.ortho}
Assume $a, b, m$ satisfy $(A1)$-$(A3)$.
If $e_1 \, , e_2$ are eigenvectors corresponding to eigenpoints $(\lambda_1, \mu_1)$ and $(\lambda_2, \mu_2)$ in $\sigma(a, b, m)$, then
\beq\label{e.ortho1}
(\lambda_2 - \lambda_1)\, b(e_1, e_2) \; + \; (\mu_2 - \mu_1)\, m(e_1 , e_2) \; = \; 0.
\eeq
\etm

\bpf
The eigenpoint equation \eqref{e.eigenequation} shows that for all $v\in V$ the following relations hold:
\[
a(e_1, v) \; = \; \lambda_1 \, b(e_1, v) + \mu_1\, m(e_1, v)\quad \text{ and }\quad
a(e_2, v) \; = \; \lambda_2 \, b(e_2, v) + \mu_2\, m(e_2, v).
\]
Take $v=e_2$ in the first equation and $v=e_1$ in the second and subtract to obtain \eqref{e.ortho1}.
\epf

When distinct eigenpoints lie on the same horizontal line, \eqref{e.ortho1} reduces to the following.

\bcor\label{c.orthoH}
Assume $a, b, m$ satisfy $(A1)$-$(A3)$, and let $\mu \in \R$ be fixed.
If $e_1, e_2$ are eigenvectors corresponding to distinct eigenpoints $(\lambda_1 , \mu)$ and $(\lambda_2 , \mu)$ in $\sigma(a, b, m)$, then
\beq\label{e.ortho2}
b(e_1, e_2)\; = \; 0 \qquad \text{and}\qquad a(e_1, e_2) \; - \; \mu \, m(e_1, e_2) \; = \; 0.
\eeq
\ecor

\bpf	Take $\mu_1 = \mu_2=\mu$ in \eqref{e.ortho1} to get $(\lambda_2 - \lambda_1)\, b(e_1, e_2) = 0$.  
Then $\lambda_2 \neq \lambda_1$ implies $b(e_1, e_2)=0$.  
Substituting this into the eigenpoint equation \eqref{e.eigenequation} gives \eqref{e.ortho2}.  
\epf

Let $E_{\lambda \mu}$ be the eigenspace in $V$ associated with the eigenpoint $(\lambda, \mu)$ of $(a, b, m)$.  
That is, let $E_{\lambda\mu}$ be the subspace generated by eigenvectors $e_{\lambda n}$ in Theorem \ref{t.abm.basis} associated to the eigenvalues $\mu_n(\lambda)$ satisfying $\mu_n(\lambda)=\mu$ for fixed $\lambda$.  
The orthogonality relations of Corollary \ref{c.orthoH} may then be expressed as
\beq\label{e.orthoH}
E_{\lambda_1\mu}\, \oplus_b\, E_{\lambda_2\mu} \qquad \text{and}\qquad E_{\lambda_1\mu}\, \oplus_{a_{\mu m}}\, E_{\lambda_2\mu}
\eeq
where $\oplus_b$ and $\oplus_{a_{\mu m}}$ indicate orthogonality with respect to $b$ and $a_{\mu m} := a -\mu m$, respectively.

In the case that distinct eigenpoints lie on the same vertical line, \eqref{e.ortho1} reduces instead to the following orthogonality relations.

\bcor\label{c.ortho3}
Assume $a, b, m$ satisfy $(A1)$-$(A3)$, and let $\lambda \in \R$ be fixed.
If $e_1, e_2$ are eigenvectors corresponding to distinct eigenpoints $(\lambda, \mu_1)$ and $(\lambda, \mu_2)$ in $\sigma(a, b, m)$, then
\beq\label{e.ortho3}
m(e_1, e_2)\; = \; 0 \qquad \text{and}\qquad a(e_1, e_2) \; - \; \lambda\, b(e_1, e_2) \; = \; 0.
\eeq
\ecor

\bpf	Take $\lambda_1 = \lambda_2=\lambda$ in \eqref{e.ortho1} to get $(\mu_2 - \mu_1)\, m(e_1, e_2) = 0$.  
Then $\mu_2 \neq \mu_1$ implies $m(e_1, e_2)=0$.  
Substituting this into the eigenpoint equation \eqref{e.eigenequation} gives \eqref{e.ortho3}.  
\epf

These results may be expressed as
\beq\label{e.orthoV}
E_{\lambda \mu_1}\oplus_m E_{\lambda \mu_2} \qquad \text{and}\qquad E_{\lambda \mu_1}\oplus_{a_{\lambda b}} E_{\lambda \mu_2}
\eeq
where $\oplus_m$ and $\oplus_{a_{\lambda b}}$ indicate orthogonality with respect to $m$ and $a_{\lambda b}$ , respectively.

From Theorem \ref{t.abm.basis}, each eigenspace $E_{\lambda\mu}$ is finite dimensional for any eigenpoint $(\lambda, \mu)$ of $(a, b, m)$, so that each linear subspace in \eqref{e.orthoH} and \eqref{e.orthoV} is finite dimensional.
These finite dimensional subspaces will be particularly useful in the following sections.


\section{Continuity of eigencurves}\label{s.continuity}
This section describes the continuity of eigencurves associated to $(a, b, m)$ and provides a first result on differentiability of eigencurves when simplicity of the eigenpoint in question is assumed. 
The following is the main result on continuity.

\btm
Assume $a, b, m$ satisfy $(A1)$-$(A3)$.
Then each variational eigencurve $\mu_n(\lambda)$ associated to $(a, b, m)$ is Lipschitz continuous.
\etm

The proof of this will be a straightforward consequence of the following result that is based on the general orthogonality result of Theorem \ref{t.ortho} of the previous section. 
Here, the open ball in $\R^2$ of radius $\varepsilon>0$ centered at $(\lambda_*, \mu_*)$ will be denoted by $B_{\varepsilon}(\lambda_*, \mu_*)$.

\btm\label{t.continuity}
Assume $a, b, m$ satisfy $(A1)$-$(A3)$, and let $(\lambda_*, \mu_*)$ be a fixed eigenpoint in $\sigma(a, b, m)$.
Then there is an $\varepsilon >0$ and a constant $C>0$ such that
\beq
\lvert \mu - \mu_*\rvert \; \leq \; C\, \lvert \lambda - \lambda_*\rvert \qquad  \text{for all }\quad (\lambda, \mu)\in B_{\varepsilon}(\lambda_*, \mu_*)\cap \sigma(a, b, m).
\eeq
\etm

\bpf
Let $\{(\lambda_i, \mu_i) : i\in \N\}$ be a sequence of eigenpoints converging to $(\lambda_*, \mu_*)$ in $\R^2$, and let $e_i$ be associated $a$-normalized eigenvectors.  
Going to a subsequence if necessary, suppose $e_i$ converges weakly to some $e_*$ in $V$.
Take $u=e_i$ in the eigenpoint equation \eqref{e.eigenequation} and let $i\rightarrow \infty$ to get
\beq\label{e.lip2}
a(e_*, v) \; = \;\lambda_*\, b(e_*, v) \; + \; \mu_*\, m(e_*, v)  \qquad \text{for each fixed }v\in V,
\eeq
using the continuity of the linear functionals associated with $a, b$, and $m$ when $v$ is fixed.  
Now let $u=v=e_i$ in \eqref{e.eigenequation} and take the limit to get 
\[
1 \; = \;  \lambda_*\, \Bc(e_*) \; + \; \mu_* \, \Mcl(e_*) \; = \; \Ac(e_*)
\]
using in this case the weak continuity of $b, m$ and taking $v=e_*$ in \eqref{e.lip2}.
This shows $e_*$ is nonzero and thus an eigenvector of $(a, b, m)$ corresponding to $(\lambda_*, \mu_*)$ as $(A1)$ and \eqref{e.lip2} hold.

Take $u=e_i$ and $v=e_i - e_*$ in \eqref{e.eigenequation} to get $\lim_{i\rightarrow \infty} [a(e_i\, , e_i -  e_*)]=0$ using the weak continuity of $b, m$.  This becomes $\lim_{i\rightarrow \infty}\Ac(e_i) \; = \; \Ac(e_*)$ upon a rearrangement which shows that $e_i$ actually converges strongly to $e_*$ in $V$ as $(A1)$ holds.

Using the orthogonality relation \eqref{e.ortho1} with $(\lambda_*, \mu_*)$ and $(\lambda_i, \mu_i)$ in place of $(\lambda_1, \mu_1)$ and $(\lambda_2, \mu_2)$ there, respectively, gives $(\lambda_i - \lambda_*)\, b(e_i\, , e_*) \; + \; (\mu_i - \mu_*)\, m(e_i \, , e_*) \; = \; 0$.
Since $\Mcl(e_*) >0$, a rearrangement of this relation and taking limits then yields
\beq\label{e.mulimit}
\lim_{i\rightarrow \infty}\frac{\mu_i - \mu_*}{\lambda_i - \lambda_*} \; = \; -\, \frac{\Bc(e_*)}{\Mcl(e_*)}.
\eeq
This implies the original sequence $\{(\lambda_i, \mu_i) :i\in \N\}$ satisfies 
\beq\label{e.diffbound}
\min_{\substack{v\in E_{\lambda_* \mu_*} \\ v\neq 0 \; \;}} \frac{-\Bc(v)}{\Mcl(v)} \quad \leq \quad
\liminf_{i\longrightarrow \infty} \, \frac{\mu_i - \mu_*}{\lambda_i - \lambda_*} \quad \leq \quad
\limsup_{i\longrightarrow \infty} \, \frac{\mu_i - \mu_*}{\lambda_i - \lambda_*} \quad \leq \quad
\max_{\substack{v\in E_{\lambda_* \mu_*} \\ v\neq 0 \; \; }} \frac{-\Bc(v)}{\Mcl(v)} 
\eeq
with the extreme sides of this relation being finite as the eigenspace $E_{\lambda_* \mu_*}$ is finite dimensional.
The desired result then follows from these estimates.
\epf

When simplicity of $\mu_n(\lambda_*)$ is assumed the following is the basic differentiability result.

\bcor
Assume $a, b, m$ satisfy $(A1)$-$(A3)$, and let $(\lambda_*, \mu_*)$ be a simple eigenpoint of $(a, b, m)$ with corresponding eigenvector $e_*$.
If $\mu_n(\lambda_*) = \mu_*$, then the eigencurve $\mu_n(\lambda)$ is differentiable at $\lambda_*$ with
\beq
\left.\frac{d\mu_n}{d\lambda}\right\vert_{\lambda_*} \; = \; -\, \frac{\Bc(e_*)}{\Mcl(e_*)}.
\eeq
\ecor
\bpf
Take $\mu_i =\mu_n(\lambda_i)$ in the previous proof.
Since $E_{\lambda_*\mu_*} = \{c\,e_* :c\in \R\}$, the inequalities in \eqref{e.diffbound} are equalities so the assertions follow as this result holds for any sequence $\lambda_i\rightarrow \lambda_*$.
\epf

Differentiability at points where the eigencurves intersect, i.e. at points where $\mu_n(\lambda_*)$ has multiplicity two or more, is more complicated, and so we devote the next section to their investigation.


\section{Differentiability results for eigencurves}\label{s.differentiability}

This section provides a local characterization of eigencurves in the case of an eigenpoint with finite multiplicity greater than one.
Specifically, eigencurves are shown to possess well-defined one-sided derivatives.
This result is then interpreted as saying that at any eigenpoint of $(a, b, m)$, the spectrum $\sigma(a, b, m)$ locally is entirely composed of smooth curves crossing at that particular point;
see Figure \ref{matrix} corresponding to the matrix problem \eqref{e.matrixproblem}, \eqref{e.matrixvalues}.

Let $(\lambda_*, \mu_*)$ be an eigenpoint of multiplicity $k_0\geq 2$ for the triple $(a, b, m)$ and let $\mu_{n_0+1}(\lambda)\; \leq \; \mu_{n_0+2}(\lambda)\;  \leq\; \cdots \; \leq \; \mu_{n_0+k_0}(\lambda)$ be the eigencurves intersecting at $(\lambda_*, \mu_*)$ so that
\[
\mu_{n_0}(\lambda_*) \; < \; \mu_* = \mu_{n_0+1}(\lambda_*) = \cdots = \mu_{n_0+k_0}(\lambda_*) \; < \; \mu_{n_0+k_0+1}(\lambda_*).
\]
Diagonalize $-b$ with respect to $m$ on $E_{\lambda_*\mu_*}$  to obtain $m$-orthonormal vectors $u_{*1}, \ldots, u_{*k_0}$ and values
\beq\label{e.diag}
b_k \; := \; \min_{Z_k}\max_{\substack{v\in Z_k \\ v\neq 0}}\frac{-\Bc(v)}{\Mcl(v)},
\eeq
with the minimization taken over all subspaces $Z_k$ of $E_{\lambda_*\mu_*}$ of dimension $k$, satisfying
\beq
b_1\leq b_2\leq \cdots \leq b_{k_0} \quad \text{ and }\quad 
b(u_{*k}\, , v) \; = \; -b_k\, m(u_{*k}\, , v) \quad \text{for all }v\in E_{\lambda_*\mu_*}
\eeq
for each $1\leq k\leq  k_0$. 
The main result on differentiability of eigencurves is the following.

\btm\label{t.diff}
Assume $a, b, m$ satisfy $(A1)$-$(A3)$.  Let $(\lambda_*, \mu_*)$ be a fixed eigenpoint in $\sigma(a, b, m)$ and define $b_k$ by \eqref{e.diag}.  
Then each eigencurve $\mu_{n_0+k}(\lambda)$ intersecting $(\lambda_*, \mu_*)$ satisfies
\beq\label{e.diff}
\lim_{\lambda \rightarrow \lambda_*^-} \frac{\mu_{n_0+k}(\lambda) - \mu_*}{\lambda - \lambda_*} \;  = \; b_{k_0-k+1}\quad\qquad\text{and}\qquad \quad
\lim_{\lambda \rightarrow \lambda_*^+} \frac{\mu_{n_0+k}(\lambda) - \mu_*}{\lambda - \lambda_*} \;  = \; b_k.
\eeq
\etm

Only the second limit in \eqref{e.diff} will be considered here as the proofs to establish each are similar.
The equality will be a direct consequence of the two associated limit extrema inequalities established below.  
For each $1\leq k\leq k_0$, define $V_{n_0+k}$ to be the subspace of $V$ generated by the vectors $e_{*1}, \ldots, e_{*n_0}, u_{*1}, \ldots, u_{*k}$, 
where $e_{*j}$ is an eigenvector associated to $\mu_j(\lambda_*)$ with $1\leq j\leq n_0$.
These {\it test subspaces based at} $\lambda_*$ will be used to prove the first of these inequalities.

\blem\label{t.limsup}
Assume $a, b, m$ satisfy $(A1)$-$(A3)$.
Let $(\lambda_*, \mu_*)$ be a fixed eigenpoint in $\sigma(a, b, m)$ and define $b_k$ by \eqref{e.diag}.
Then each eigencurve $\mu_{n_0+k}(\lambda)$ intersecting $(\lambda_*, \mu_*)$ satisfies
\beq
\limsup_{\lambda\,\rightarrow \,\lambda_*^+}\, \frac{\mu_{n_0+k}(\lambda) - \mu_*}{\lambda - \lambda_*}\;\; \leq \; \; b_k.
\eeq
\elem

\bpf
Let $\mu_{n_0+k}(\lambda)$ be a fixed eigencurve intersecting $(\lambda_*, \mu_*)$ and assume the contrary.  
That is, let $\varepsilon >0$ and assume there are values $\lambda_j >  \lambda_*$ decreasing to $\lambda_*$ satisfying 
\[
\mu_{n_0 + k}(\lambda_j)  \; \; \geq \;  \mu_* \; + \; (b_k+ \varepsilon)(\lambda_j - \lambda_*).
\]
Consider $V_{n_0+k}$ as one of the $(n_0+k)$-dimensional subspaces in the variational characterization \eqref{e.minmax} for $\mu_{n_0 +k}(\lambda_j)$ for each $j\in \N$ to obtain $m$-normalized vectors $v_j\in V_{n_0+k}$ satisfying
\beq\label{e.contra1}
\Ac(v_j) \; - \; \lambda_j\, \Bc(v_j) \; \; \geq \; \; \mu_* \; + \; (b_k+\varepsilon)(\lambda_j - \lambda_*)
\eeq
upon using the continuity of $\Ac - \lambda_j\Bc$ on the $m$-unit sphere of $V_{n_0+k}$.
Extract a subsequence of the $v_j$, denoting it again by $v_j$, that converges strongly to a vector $v_*$ in $V_{n_0+k}$.
This vector satisfies $\Mcl(v_*)=1$, and also $\Bc(v_j)\rightarrow \Bc(v_*)$ and $\Ac(v_j)\rightarrow \Ac(v_*)$.
Taking the limit in \eqref{e.contra1} shows that $v_*$ also satisfies $\Ac(v_*) \; - \; \lambda_*\, \Bc(v_*) \; \; \geq \; \; \mu_*$.
From \eqref{e.minmax}, the inequality
\beq\label{e.contra2}
\Ac(v) \; - \; \lambda_*\, \Bc(v) \; \; \leq \; \;  \mu_* \, \Mcl(v) 
\eeq
holds for all $v\in V_{n_0+k}$ with equality being satisfied for all $v\in E_{\lambda_*\mu_*}\cap V_{n_0+k}= \text{span}\{u_{*1}, \ldots, u_{*k}\}$.
These last two estimates imply $v_*$ belongs to $\text{span}\{u_{*1}, \ldots, u_{*k}\}$.

Rewrite inequality \eqref{e.contra1} as
\[
\Ac(v_j) \; - \; \lambda_*\, \Bc(v_j) \; + \; (\lambda_* - \lambda_j)\, \Bc(v_j) \; \; \geq \; \; \mu_* \; + \; (b_k +\varepsilon)(\lambda_j - \lambda_*).
\]
Use the inequality \eqref{e.contra2} for each $v_j$ to simplify the above line to $-\Bc(v_j) \geq b_k\, +\varepsilon$ using the $m$-normalization of the $v_j$.  Then get $-\Bc(v_*)  \geq  b_k+\varepsilon$ in the limit.
This, however, contradicts the diagonalization \eqref{e.diag} of $-b$ with respect to $m$ on $E_{\lambda_*\mu_*}$.

Therefore, there is a $\delta>0$ such that for $0 < \lambda - \lambda_* < \delta$,
\[
\max_{\substack{v\in V_{n_0+k} \\ v\neq0\;\;}} \left\{  \frac{\Ac(v) - \lambda\, \Bc(v)}{\Mcl(v)}\right\} \; \; \leq \; \; \mu_* \; + \; (b_k + \varepsilon)(\lambda - \lambda_*).
\]
Using \eqref{e.minmax} for $\mu_{n_0 +k}(\lambda)$ then gives $\mu_{n_0+k}(\lambda)\, \leq  \, \mu_* +(b_k +\varepsilon)(\lambda - \lambda_*)$ which leads to
\[
 \limsup_{\lambda\, \rightarrow \, \lambda_*^+}\left( \frac{\mu_{n_0+k}(\lambda) - \mu_*}{\lambda - \lambda_*}\right) \; \; \leq \; \; b_k+\varepsilon
\]
upon a rearrangement.
Since $\varepsilon >0$ was arbitrary the desired result then holds.
\epf

To establish the second inequality associated with the second limit in \eqref{e.diff}, the following lemma will be used to extract specific convergent subsequences from `neighboring' eigenspaces of $E_{\lambda_*\mu_*}$.

\blem\label{l.extract}
Assume $a, b, m$ satisfy $(A1)$-$(A3)$, and let $\mu_n(\lambda)$ be a fixed eigencurve intersecting the eigenpoint $(\lambda_*, \mu_*)$.
Suppose $\lambda_i \rightarrow \lambda_*$ and let $\{e_{in} : i\in \N\}$ be a sequence of $m$-normalized eigenvectors associated to $\mu_n(\lambda_i)$.
Then there is a subsequence of $\{e_{in} : i\in \N\}$  that converges strongly to an eigenvector $v_{*n}$ associated to $\mu_n(\lambda_*)$ satisfying $\Mcl(v_{*n}) =1$.
\elem

\bpf
Let $r_0 > \sup\{\lvert \lambda_i\rvert : i\in \N\}$ to obtain a $\tau>0$ as in Theorem \ref{t.tau}.
The bilinear forms $a_{\lambda_i, \tau} : = a - \lambda_i\, b + \tau\, m$ satisfy the identity $a_{\lambda_i, \tau}(e_{in}\, , e_{in}) \; = \; \mu_n(\lambda_i) \; + \; \tau$ upon using the eigenpoint equation \eqref{e.eigenequation} with $u=e_{in}$.
Equation \eqref{e.alamtau} then gives
\[
\mu_n(\lambda_i) + \tau \; \geq \; \frac{1}{2}\Ac(e_{in}) 
\] 
for each $i\in\N$.
This and continuity of $\mu_n(\lambda)$ on $[-r_0, r_0]$ imply $\{e_{in}: i\in \N\}$ is a bounded set and thus a weakly precompact set in $V$ as $(A1)$ holds.
Extract a subsequence, denoted again by $\{e_{in} : i\in \N\}$, weakly converging to a vector $v_{*n}$ in $V$.
Since $\Mcl(e_{in}) = 1$ for each $i\in \N$, it follows that $\Mcl(v_{*n}) =1$.
Take $i\rightarrow \infty$ in the eigenpoint equation \eqref{e.eigenequation} with $u=e_{in}$ to see that $v_{*n}$ is an eigenvector associated to $\mu_n(\lambda_*)$.

Since $\|v_{*n} - e_{in}\|_a^2 = \Ac(v_{*n}) - 2\, a(v_{*n}\, , e_{in}) + \Ac(e_{in})$, it suffices to show $\Ac(e_{in})\rightarrow \Ac(v_{*n})$ to establish strong convergence of $e_{in}$ to $v_{*n}$.
Therefore, from
\[
\Ac(v_{*n}) \; = \; \lambda_*\, \Bc(v_{*n}) + \mu_n(\lambda_*)\, \Mcl(v_{*n}) \; = \; \lim_{i\rightarrow\infty}\big( \lambda_i\, \Bc(e_{in}) + \mu_n(\lambda_i)\, \Mcl(e_{in})\big) \; = \; \lim_{i\rightarrow \infty}\Ac(e_{in})
\]
the assertions of the lemma hold.
\epf

\blem
Assume $a, b, m$, satisfy $(A1)$-$(A3)$.  
Let $(\lambda_*, \mu_*)$ be a fixed eigenpoint in $\sigma(a, b, m)$ and define $b_k$ by \eqref{e.diag}.
Then each eigencurve $\mu_{n_0+k}(\lambda)$ intersecting $(\lambda_*, \mu_*)$ satisfies
\beq
\liminf_{\lambda \, \rightarrow \, \lambda_*^+} \, \frac{\mu_{n_0+k}(\lambda) - \mu_*}{\lambda - \lambda_*} \; \; \geq \; \; b_k.
\eeq
\elem

\bpf
Let $\mu_{n_0+k}(\lambda)$ be a fixed eigencurve intersecting $(\lambda_*, \mu_*)$.
Let $\varepsilon>0$ and assume there are values $\lambda_i> \lambda_*$ decreasing to $\lambda_*$ satisfying
\beq\label{e.contra3}
\mu_{n_0+k}(\lambda_i) \; \; \leq \; \; \mu_* \; + \; (b_k - \varepsilon)(\lambda_i - \lambda_*).
\eeq
Using the eigencurves $\mu_n(\lambda)$ with $n_0+1\leq n\leq  n_0+k$ we show \eqref{e.contra3} leads to a contradiction.
For such index values $n$, let $e_{in}$ be an eigenvector associated to $\mu_n(\lambda_i)$, taken here to be $m$-normalized.
Use Lemma \ref{l.extract} to extract subsequences from the union of all such $e_{in}$, denoting them simply by $\{e_{in} : i\in \N\}$, and $m$-orthonormal vectors $v_{*n}$ associated to $\mu_n(\lambda_*)$, such that $e_{in}$ converges strongly to $v_{*n}$ in $V$.
By construction, $\text{span}\{v_{*n} : n_0+1\leq n\leq n_0+k\}$ is equal to $\text{span}\{u_{*1}, u_{*2}, \ldots , u_{*k}\}$. 
Using the orthogonality relation of Theorem \ref{t.ortho} for each $n_0+1\leq n\leq n_0+k$ gives
 \[
 \lim_{i\rightarrow \infty} \, \frac{\mu_n(\lambda_i) - \mu_n(\lambda_*)}{\lambda_i - \lambda_*} \; \; = \; \; \lim_{i\rightarrow \infty}\, \frac{-b(e_{in}, v_{*n})}{m(e_{in}, v_{*n})} \; \; = \; \; -\Bc(v_{*n}).
 \]
This and \eqref{e.contra3} show  $-\Bc(v_{*n_0+k}) \leq  b_k - \varepsilon$, so the ordering $b_1\leq b_2\leq \cdots \leq b_{k_0}$ then implies
\[
\frac{-\Bc(v)}{\Mcl(v)} \; \; \leq \; \; b_k - \varepsilon
\]
for all nonzero $v$ in $\text{span}\{v_{*n} : n_0+1\leq n\leq n_0+k\}$, contradicting the variational characterization \eqref{e.diag} of $b_k$.  
The result then follows similarly as in the proof of Theorem \ref{t.limsup}.
\epf

The result of Theorem \ref{t.diff} says that an eigencurve $\mu_{n_0+k}(\lambda)$ intersecting an eigenpoint $(\lambda_*, \mu_*)$ of multiplicity $k_0$ greater than one, possesses well-defined one-sided derivatives at that eigenpoint given by \eqref{e.diff}.  
Since the left- and right-derivative values $b_{k_0-k+1}$ and $b_k$ are not equal to each other in general at such points, this implies such an eigencurve $\mu_{n_0+k}(\lambda)$ is in general not differentiable at such an eigenpoint.

However, since the eigencurve $\mu_{k_0-k+1}(\lambda)$ has $b_k$ as its left-derivative at $(\lambda_*, \mu_*)$, which is equal to the right-derivative of $\mu_{n_0+k}(\lambda)$ at this eigenpoint, the theorem implies the spectrum of $(a, b, m)$ can be characterized locally around $(\lambda_*, \mu_*)$ as being composed of $k_0$ curves that are differentiable at $(\lambda_*, \mu_*)$.  
These differentiability results are observed in the spectrum shown in Figure \ref{matrix} corresponding the specific matrix problem \eqref{e.matrixproblem}, \eqref{e.matrixvalues}.


\section{Asymptotic results}\label{s.asymptotic}

The previous two sections dealt with regularity issues pertaining to eigencurves.  
In the next few sections, the interest will now focus on providing a geometrical description of these variational eigencurves.

First, we shall treat the asymptotic behaviour of eigencurves.
In particular, it is shown that the dimension of subspaces on which the quadratic form $\Bc$ is strictly positive, or strictly negative, controls the number of eigencurves that go asymptotically down to negative infinity as $\lambda\rightarrow \infty$.
Using a more careful analysis, explicit formulae are obtained to describe the asymptotic behaviour, the results depending on spectral data for the pair $(b, m)$ of bilinear forms.

For the first eigencurve, the result is the following.

\blem
Assume $a, b, m$ satisfy $(A1)$-$(A3)$, and let $\mu_1(\lambda)$ be the first variational eigencurve associated with $(a, b, m)$.\\
$(i)$\quad If there is a vector $\hat{v}\in V$ such that $\Bc(\hat{v})>0$, then
\[
\lim_{\lambda\rightarrow \infty}\mu_1(\lambda) \; = \; -\infty.
\]
$(ii)$\quad If there is a vector $\hat{w}\in V$ such that $\Bc(\hat{w})<0$, then 
\[
\lim_{\lambda\rightarrow\, -\infty}\mu_1(\lambda) \; = \; -\infty.
\]
\elem

\bpf
The variational characterization \eqref{e.minmax} for $\mu_1(\lambda)$ gives
\[
\mu_1(\lambda) \; \; \leq \; \; \frac{\Ac(v) - \lambda\, \Bc(v)}{\Mcl(v)} \qquad \text{for all nonzero }v\in V.
\]
For such $v$ fixed, the right side of this inequality is an affine function in $\lambda$ with the sign of its slope determined by the sign of $\Bc(v)$.
Considering $\hat{v},\hat{w}$ as in the theorem and taking the appropriate limit in $\lambda$ gives the desired results.
\epf

This result may be interpreted as saying that the sign of the quadratic form $\Bc$ on one-dimensional subspaces (spanned by such $\hat{v}$ or $\hat{w}$) may be used to determine the asymptotic behaviour of the first variational eigencurve $\mu_1(\lambda)$.  
The next result says that the asymptotic behaviour of the $k^{th}$ eigencurve is determined by the sign of $\Bc$ on $k$-dimensional subspaces.

\blem
Assume $a, b, m$ satisfy $(A1)$-$(A3)$, and let $\mu_k(\lambda)$ be the $k^{th}$-variational eigencurve associated with $(a, b, m)$.\\
$(i)$\quad If $\hat{v}_1, \ldots, \hat{v}_k$ are $k$ linearly independent vectors at which $\Bc$ is strictly positive, then
\[
\lim_{\lambda\rightarrow \infty} \mu_k(\lambda) \; = \; -\infty.
\]
$(ii)$\quad If $\hat{w}_1, \ldots, \hat{w}_k$ are $k$ linearly independent vectors at which $\Bc$ is strictly negative, then 
\[
\lim_{\lambda\rightarrow \, - \infty}\mu_k(\lambda) \; = \; - \infty.
\]
\elem

\bpf
Let $\lambda>0$ and let $\hat{Z}_k$ denote the subspace spanned by the vectors $\hat{v}_1, \ldots, \hat{v}_k$ at which $\Bc$ is stricly positive.  The variational characterization \eqref{e.minmax} for $\mu_k(\lambda)$ gives
\[
\mu_k(\lambda) \; \; \leq \; \; \sup_{\substack{v\in \hat{Z}_k \\ v\neq 0}}\frac{\Ac(v) - \lambda\, \Bc(v)}{\Mcl(v)} \; \; \leq \; \; 
\sup_{\substack{v\in \hat{Z}_k \\v\neq 0}} \left\{ \frac{\Ac(v)}{\Mcl(v)}  \right\}\; + \; \lambda\, \sup_{\substack{v\in \hat{Z}_k \\ v\neq 0}}\left\{ - \frac{\Bc(v)}{\Mcl(v)}\right\}.
\]
This holds as $\Bc$ is strictly positive on the $m$-unit sphere in the finite dimensional space $\hat{Z}_k$.
The right side of this relation is an affine function in $\lambda$ .
Part $(i)$ then follows upon taking the limit, and part $(ii)$ is similarly proved.
\epf

When $\Bc$ is strictly positive on some $k$-dimensional subspace of $V$, it follows from this result that the first $k$ eigencurves asymptotically decrease to $-\infty$ as $\lambda$ increases to $+\infty$.  
The positivity of $\Bc$ is easily verifiable, so this asymptotic information for eigencurves is straightforward to establish.
A similar statement can be made when $\Bc$ is strictly negative on a finite dimensional subspace.

To determine a more precise asymptotic description, define the values
\beq\label{e.bmspectral}
\eta_k\; := \; \inf_{Z_k}\sup_{\substack{v\in Z_k \\ v\neq 0}} \left\{\, - \frac{\Bc(v)}{\Mcl(v)} \right\} 
\eeq
where the infimum is taken over all $k$-dimensional subspaces $Z_k$ of $V$.
It is worth noting that these values may equal $-\infty$, but never $+\infty$.
The following theorem shows that the asymptotic behaviour of eigencurves are prescribed by these spectral values $\eta_k$, and that the sign of $\Bc$ plays a role here as in the previous lemmas.

\btm
Assume $a, b, m$ satisfy $(A1)$-$(A3)$.
Then the $k^{th}$-eigencurve $\mu_k(\lambda)$ satisfies
\beq\label{e.asymptotek}
\lim_{\lambda\rightarrow \infty}\frac{\mu_k(\lambda)}{\lambda} \; \; = \; \; \eta_k.
\eeq
\etm

\bpf
Assume $\lambda>0$ and let $Z_k$ be a $k$-dimensional subspace of $V$.
Then
\[
\sup_{\substack{v\in Z_k \\ v\neq 0}}\, \frac{\Ac(v) -  \lambda\, \Bc(v)}{\Mcl(v)} \; \; \geq \; \; \lambda\, \sup_{\substack{v\in Z_k \\ v\neq 0}} \left\{ - \frac{\Bc(v)}{\Mcl(v)}\right\}.
\]
Taking the infimum over all such $Z_k$, and taking the limit inferior after rearranging $\lambda$ gives 
\[
\liminf_{\lambda\rightarrow \infty} \, \frac{\mu_k(\lambda)}{\lambda} \; \; \geq \; \; \eta_k,
\]
using the variational characterization of $\mu_k(\lambda)$. 

To establish the desired reverse inequality, consider first the case in which $\eta_k >-\infty$.
Let $\varepsilon >0$ and take $\hat{Z}_k$ to be a $k$-dimensional subspace of $V$ satisfying
\[
\sup_{\substack{v\in \hat{Z}_k \\  v\neq 0}}\left\{ - \frac{\Bc(v)}{\Mcl(v)}\right\} \; \; \leq \; \; \eta_k \; + \; \varepsilon.
\]
The variational characterization \eqref{e.minmax} for $\mu_k(\lambda)$ shows that
\[
\frac{\mu_k(\lambda)}{\lambda} \; \; \leq \; \; \frac{1}{\lambda}\sup_{\substack{v\in \hat{Z}_k \\ v\neq 0}}\left\{ \frac{\Ac(v)}{\Mcl(v)} \right\} \; + \; \eta_k \; + \; \varepsilon
\]
which gives $\limsup_{\lambda \rightarrow \infty} \mu_k(\lambda)/\lambda \;  \leq  \; \eta_k \; + \; \varepsilon$.
Taking $\varepsilon$ down to zero then gives the desired converse inequality so that \eqref{e.asymptotek} holds in the case $\eta_k$ is finite. 

When $\eta_k=-\infty$, for $n\in \N$ take $Z_{k, n}$ to be a $k$-dimensional subspace of $V$ such that
\[
\sup_{\substack{v\in Z_{k, n} \\  v\neq 0}}\left\{ - \frac{\Bc(v)}{\Mcl(v)}\right\} \; \; \leq \; \;  -n.
\]
The characterization \eqref{e.minmax} for the eigencurve $\mu_k(\lambda)$ gives in this case
\[
\frac{\mu_k(\lambda)}{\lambda} \; \; \leq \; \; \frac{1}{\lambda}\sup_{\substack{v\in Z_{k, n} \\ v\neq 0}}\left\{ \frac{\Ac(v)}{\Mcl(v)} \right\} \; - \; n.
\]
Taking the limit superior in this inequality and then taking $n$ to $+\infty$ gives the desired result.
\epf

For the specific matrix problem \eqref{e.matrixproblem},\eqref{e.matrixvalues}, it is easy to see that the quadratic form $\Bc$ corresponding to the matrix $B$ is striclty positive on a two-dimensional subspace of $\mathbb{R}^3$.  
Hence, as seen in Figure \ref{matrix}, the first and second eigencurves asymptotically go to $-\infty$ as $\lambda\rightarrow \infty$;
the eigencurves are asymptotic to lines whose slopes are the negative eigenvalues of $B$.

Richardson's equation is discussed in Section 3 of Binding and Vokmer \cite{BV96} to exemplify their results on Sturm-Liouville problems (equations \eqref{e.SL},\eqref{e.SLBC} here).
The quadratic form $\Bc$ associated with this example is $\Bc(y) \; := \int_{-1}^1 \text{sgn}(x)\, y(x)^2\, dx$ defined on $L^2[-1, 1]$.
It is obvious that this form is strictly positive on an infinite dimensional subspace of $L^2[-1, 1]$ so that all of the eigencurves $\mu_n(\lambda)$ for this example decrease to $-\infty$ as $\lambda\rightarrow +\infty$.
Likewise, $\mu_n(\lambda)\rightarrow -\infty$ as $\lambda\rightarrow -\infty$ for this example.
This behaviour is seen in Figure 3.1 of \cite{BV96}.
In general, their Theorem 2.2 provides explicit formulae for the asymptotic behaviour of the eigencurves in terms of the extreme values the weight function $r=r(x)$ takes on the interval $[t_0, t_1]$ (appearing in equation \eqref{e.SL} of Section \ref{s.examples} here).

Recall that for the Robin-Steklov problem \eqref{e.RS}, the function $b_0$ may be sign-changing.
If $b_0$ is strictly positive $\sigma$-$a.e.$ on a subset of $\partial\Omega$ of strictly positive Hausdorff measure, then all of the eigencurves $\mu_n(\lambda)$ for this example decrease to $-\infty$ as $\lambda\rightarrow \infty$.
If there is a subset of $\partial\Omega$ of strictly positive Hausdorff measure where $b_0$ is strictly negative $\sigma$-$a.e.$, then in this case all of the eigencurves decrease to $-\infty$ as $\lambda\rightarrow -\infty$.
In contrast to the Sturm-Liouville problem \eqref{e.SL},\eqref{e.SLBC} where the $b$ form is a weighted version of the $m$ bilinear form there (the $L^2$-inner product), the $b, m$ bilinear forms for the Robin-Steklov problem are not variants of each other - one comprises a boundary integral whereas the other an interior-region integral.  In this case, the explicit asymptotic behaviour of Robin-Steklov eigencurves is described in terms of spectral data $\eta_k$ given by \eqref{e.bmspectral} where the pair $(b, m)$ is given by \eqref{e.bmforms}; $\eta_k$ is infinite in this case.


\vspace{1em}
\section{Straight lines within the spectrum and linear independence}\label{s.independence}

In this section, the degeneracy of the form $b$ is used to determine when straight lines appear within the spectrum of $(a, b, m)$.  
When this holds, a further orthogonal decomposition is found for some of the associated eigenspaces.
Related to these results is a description of linear independence of vectors associated with eigenspaces corresponding to eigenpoints lying on the same horizontal line.

To proceed with the analysis, recall that the {\it null space} $N(b)$ of a bilinear form $b$ obeying (A2) is the set of all vectors $u\in V$ satisfying
\begin{equation}
b(u, v) \; = \; 0 \qquad \text{for all }v\in V.
\end{equation}
When $N(b) =\{0\}$ then $b$ is said to be {\it non-degenerate}, otherwise $b$ is said to be {\it degenerate}.  
It is worth noting that the degeneracy of $b$ has played no role in the preceding analysis.

\btm\label{t.independent}
Assume $a, b, m$ satisfy $(A1)$-$(A3)$.
Let $\mu_*\notin \sigma(a, m)$ and $(\lambda_i, \mu_*)$ be distinct eigenpoints in $\sigma(a, b, m)$ for $1\leq i\leq k$, with $k\geq 2$.
If $e_i$ are eigenvectors associated with $(\lambda_i, \mu_*)$, then the set $\Ec_k := \{e_1, \ldots,e_k\}$ is linearly independent.
\etm

\bpf
Suppose $\Ec_k$ is linearly dependent. 
Then without loss of generality we may assume $e_k \in \text{span}\{e_1, \ldots, e_{k-1}\}$.
If $e_k=\alpha\, e_i $ for some $1\leq i\leq k-1$ and some nonzero constant $\alpha$, then by subtracting the eigenpoint equations that $e_k, e_i$ satisfy shows that $b(e_k, v) =0$ for all $v\in V$.
This says $e_k\in N(b)$ which implies $(\lambda, \mu_*)\in \sigma(a, b, m)$ for any $\lambda\in \R$.
In particular, $(0, \mu_*)\in \sigma(a, b, m)$ or that $\mu_*\in \sigma(a, m)$, in contradiction with $\mu_*\notin \sigma(a, m)$.

Without loss of generality assume now that $e_k = \sum_{i=1}^{k-1}\alpha_i\, e_i$ with $\alpha_i\neq 0$ for all $i$ and that the dimension $d$ of the span of $\Ec_k$ satisfies $2\leq d= k-1$.
Substituting this expression for $e_k$ into the eigenpoint equation that $e_k$ satisfies gives
\[
\sum_{i=1}^{k-1}\alpha_i\, a(e_i, v) \; = \; \lambda_k\, b(e_k, v) \; + \; \mu_*\, \sum_{i=1}^{k-1}\alpha_i\, m(e_i, v)
\]
holding for all $v\in V$.
Considering the terms that involve the $a, m$ bilinear forms, and using the eigenpoint equations for $e_1, \ldots, e_{k-1}$ simplifies this expression to
\[
\sum_{i=1}^{k-1}\alpha_i\, \lambda_i\, b(e_i, v) \; = \; \lambda_k\, b(e_k, v).
\]
Using the expansion for $e_k$ and rearranging terms reduces this last relation to
\beq\label{e.tildee}
b\left(\tilde{e}\, ,\, v\right) =0 \qquad \quad \text{with }\quad \tilde{e}\; :=\;  \sum_{i=1}^{k-1}\alpha_i\, (\lambda_i - \lambda_k)\, e_i.
\eeq
The vector $\tilde{e}$ is nontrivial since the $\lambda_i$ are distinct and $\{e_1, \ldots, e_{k-1}\}$ is linearly independent.
As \eqref{e.tildee} holds for all $v\in V$, it follows that $\tilde{e}$ belongs to $N(b)$.

Now, for a given value $\lambda\in \R\setminus\{\lambda_1, \ldots, \lambda_{k-1}\}$, define the vector
\[
e_{\lambda} \; :=\; \sum_{i=1}^{k-1}\alpha_i\left( \frac{\lambda_i - \lambda_k}{\lambda_i - \lambda}\right) e_i.
\]
For fixed $v\in V$, the expression $\lambda\, b(e_{\lambda}\, , v) + \mu_*\, m(e_{\lambda}\, , v)$ is then equal to
\beq\label{e.expert}
\sum_{i=1}^{k-1}\alpha_i\left( \frac{\lambda_i - \lambda_k}{\lambda_i - \lambda}\right) \, [\lambda \, b(e_i, v) \; + \; \mu_*\, m(e_i, v)]
\eeq
upon collecting like terms.
Since the quantity $\lambda\, b(e_i, v) \, + \, \mu_*\, m(e_i, v)$ is equal to the quantity $\lambda_i\, b(e_i, v) \, + \, \mu_*\, m(e_i, v) \, + \, (\lambda - \lambda_i)b(e_i, v)$, it follows from the eigenpoint equation that $e_i$ satisfies that the expression \eqref{e.expert} is equal to
\[
\sum_{i=1}^{k-1}\alpha_i\left( \frac{\lambda_i - \lambda_k}{\lambda_i - \lambda}\right) \, [a(e_i, v) \; + \; (\lambda - \lambda_i)\, b(e_i, v)].
\]
Distributing and using the definitions of $\tilde{e}, e_{\lambda}$ shows that this last expression is equal to the quantity $a(e_{\lambda}, v) \, + \, b(\tilde{e}, v)$, which reduces to $a(e_{\lambda}, v)$ since $\tilde{e}\in N(b)$.
This shows $e_{\lambda}$ satisfies $a(e_{\lambda}, v) = \lambda\, b(e_{\lambda}, v) + \mu_*\, m(e_{\lambda}, v)$ for all all $v\in V$, so that $e_{\lambda}$ is an eigenvector of $(a, b, m)$ corresponding to the eigenvalue $(\lambda, \mu_*)$.
Taking $\lambda=0$ implies that $(0, \mu_*)\in \sigma(a, b, m)$ which contradicts $\mu_*\notin \sigma(a, m)$.
Therefore, $\Ec_k$ is a linearly independent set.
\epf

The next results considers the case $\mu_*\in \sigma(a, m)$, and it says that when there is a nontrivial intersection of the corresponding eigenspace $E_{0\mu_*}$ and the null space $N(b)$ of $b$, then there is an entire straight line within the spectrum of $(a, b, m)$.

\btm
Assume $a, b, m$ satisfy $(A1)$-$(A3)$, respectively, and let $\mu_*\in \R$ be fixed.
Then any nonzero vector  $e\in E_{0 \mu_*}\cap N(b)$ also belongs to the eigenspace $E_{\lambda\mu_*}$ for each $\lambda\in \R$, i.e., the straight line
$\{(\lambda, \mu_*) : \, \lambda\in \R\}$ is contained in $\sigma(a, b, m)$ in this case.
\etm

\bpf
Such an $e$ satisfies $a(e, v) = \mu_*\, m(e, v)$ and also $b(e, v)=0$ for all $v\in V$. 
Multiply the second relation by $\lambda$ and add this to the first relation to get that $e$ is an eigenvector of $(a, b, m)$ corresponding to $(\lambda, \mu_*)$ for any $\lambda$.  
The assertions then follow.
\epf

Let $N_{0\mu_*}(b)$ denote the subspace $E_{0\mu_*}\cap N(b)$ for fixed $\mu_*\in \R$.
We point out that in many cases $N_{0\mu_*}(b)$ is simply the trivial subspace.
When it is nontrivial, denote by $\hat{E}_{\lambda\mu_*}$ the $m$-orthogonal complement of $N_{0\mu_*}(b)$ in $E_{\lambda\mu_*}$ to obtain the decomposition
\beq
E_{\lambda\mu_*} \; = \; \hat{E}_{\lambda\mu_*}\, \oplus_m \, N_{0\mu_*}(b)
\eeq
where $\oplus_m$ indicates orthogonality with respect to the $m$ inner product on $V$.
This yields a further orthogonal decomposition of eigenspace to those given in Section \ref{s.ortho}.

It is instructive at this point to provide some easy concrete examples of spectra to illustrate the preceding results.  
In particular, we describe spectra consisting entirely of (not necessarly horizontal) straight lines.

Suppose $a, m$ satisfy (A1), (A3), and let $\Lambda_0 : = \{\mu_{0n} : n\in\mathbb{N}\}$ and $\Ec_0 :=\{e_{0n} : n\in \mathbb{N}\}$ be the eigendata described in Section \ref{s.definitions}  satisfying
\[
a(e_{0n}\, , v) \; = \; \mu_{0n} \, m(e_{0n}\, , v) \qquad \text{ for all } v\in V.
\]
Consider the form $b_{\varepsilon} := \varepsilon \, m$ defined on $V$ for fixed $\varepsilon\in \mathbb{R}$.  
When $\varepsilon =0$, the pair $(a, b_{\varepsilon})$ has no eigenvalues and then the horizontal straight lines $\mu_n(\lambda) = \mu_{0n}$ completely make up the spectrum $\sigma(a, b_{\varepsilon}, m)$.

When $\varepsilon \neq 0$,  then $\{\varepsilon^{-1}\mu_{0n} :n\in \mathbb{N}\}$ is the sequence of eigenvalues of $(a, b_{\varepsilon})$ so that the spectrum $\sigma(a, b_{\varepsilon}, m)$ comprises precisely the nonhorizontal straight lines $\mu_n(\lambda) = \mu_{0n} - \varepsilon \, \lambda$.
In these two examples, the functional forms of the eigencurves are linear.

As a third example, take $K\subset \mathbb{N}$ to be an indexing set of finite cardinality $\lvert K\rvert$, and let $\varepsilon$ be a vector in $\mathbb{R}^{\lvert K\rvert}$ with nonzero entries $\varepsilon_k$.  Consider the bilinear form $b_{\varepsilon K}$ defined by
\beq
b_{\varepsilon K}(u, v) \; := \; \sum_{k\in K} \varepsilon_k\, m(e_{0k}\, , u)\, m(e_{0k}\,, v) \qquad \text{for }u, v\in V.
\eeq
If $n\notin K$, then $b_{\varepsilon K}(e_{0n}\, , v) =0$ for all $v\in V$.  
Using the eigenpoint equation \eqref{e.eigenequation} for $(a, b_{\varepsilon K}, m)$ it is easy to see that the graph of the constant function $\lambda \mapsto \mu_{0n}$ belongs to $\sigma(a, b_{\varepsilon K}, m)$.
When $n\in K$, then the eigenpoint equation shows that the graph of the function $\lambda \mapsto \mu_{0n} - \varepsilon_n \lambda$ belongs to $\sigma(a, b_{\varepsilon K}, m)$.  
Here, the nonhorizontal lines are not necessarily parallel nor have slopes of the same sign.  
Therefore, in this example the eigencurves $\mu_n(\lambda)$ for this system are instead piecewise-linearly defined.

From these examples, it is easy to see that the degeneracy of $b$ plays a major role in the geometry of the spectrum of $(a, b, m)$ and also in the asymptotic behaviour of the eigencurves.
We also remark that since only straight lines through points $(0, \mu_*)$, with $\mu_*\in \sigma(a, m)$, can be in $\sigma(a,b,m)$, there are at most countably many straight lines in spectrum of $(a, b, m)$.


\section{Intersecting the spectrum with straight lines}\label{s.geometry}

This section deals with the geometrical problem of establishing an upper bound for the number of {\it components} the $n^{th}$-eigencurve may have above a given a horizontal straight line $\mu = \mu_*$.
The case where the intersecting line is not horizontal is discussed at the end of the section, and the results are seen to also hold for this case as a consequence of a simple change-of-variables.

For a fixed eigencurve $\mu_n(\lambda)$ and fixed $\mu_*\in \R$, denote by $\{\mu_n > \mu_*\}$ the superlevel set of $\mu_n(\lambda)$ given by $\{\lambda\in \R : \mu_n(\lambda) > \mu_*\}$.  
A finite interval $(s, t)\subset \R$ is said to be a {\it component} of $\{\mu_n > \mu_*\}$ provided $\mu_n(s) = \mu_n(t) =\mu_*$ and $\mu_n(\lambda) > \mu_*$ for all $s< \lambda < t$.  
A semi-infinite interval $(s, \infty) \subset \R$ is also said to be a {\it component} of $\{\mu_n> \mu_*\}$ provided $\mu_n(s)=\mu_*$  and $\mu_n(\lambda)>\mu_*$ for all $\lambda\in (s, \infty)$;
a component of the form $(-\infty, t)$ is similarly defined.
We shall often say a component $(s,t)$ of $\{\mu_n >\mu_*\}$ is a component of $\mu_n(\lambda)$ above $\mu_*$.

For the first eigencurve, the geometrical result is the following.

\blem\label{l.curve1}
Assume $a, b, m$ satisfy $(A1)$-$(A3)$ and let $\mu_*\notin\sigma(a,m)$ be fixed.
Then  $\{\mu_1 > \mu_*\}$ has at most one component.
\elem

\bpf
This follows from the concavity of $\mu_1(\lambda)$ given by Lemma \ref{l.mu1concave}.
\epf

To begin the geometrical description of higher eigencurves, we first show in the next lemma that the eigencurves satisfy a {\it nested property} and then in Lemma \ref{l.eprobe} we provide a {\it sign-condition} that eigenvectors associated with endpoints of a component must satisfy.

\blem\label{l.nested}
Assume $a, b, m$ satisfy $(A1)$-$(A3)$, and let $\mu_*\notin\sigma(a,m)$ and $n,k\in \mathbb{N}$ be fixed with $k>n$.
Let $I_n$ be a component of $\{\mu_n >\mu_*\}$ and suppose $\{\mu_k > \mu_*\}$ has a component.
Then for each integer $n\leq i \leq k$, there is a component $I_i$ of $\{\mu_i > \mu_*\}$ such that these components satisfy 
$I_n \subset I_{n+1} \subset \cdots \subset I_i\subset \cdots \subset I_k$.
\elem

\bpf
By definition of eigencurves, $\mu_k(\lambda)\geq \mu_n(\lambda)$ for each $\lambda\in \R$, so $\mu_k(\lambda)> \mu_*$ on $I_n$.
Since $\{\mu_k > \mu_*\}$ has at least one component, there is a component $I_k$ of $\{\mu_k>\mu_*\}$ containing $I_n$.
Since $\mu_n(\lambda) \leq \mu_i(\lambda)\leq \mu_k(\lambda)$ holds for $n\leq i \leq k$, again by definition, it follows that $\{ \mu_i > \mu_*\}$ contains a component $I_i$ satisfying $I_k \supset I_i \supset I_n$ and the assertions then follow.
\epf

It is worth noting that the components in Lemma \ref{l.nested} may be semi-infinite intervals.

\blem\label{l.eprobe}
Assume $a, b, m$ satisfy $(A1)$-$(A3)$ and let $\mu_*\notin \sigma(a,m)$ be fixed.
If the bounded interval $(s, t)$ is a component of $\{\mu_n > \mu_*\}$, then eigenvectors $e_{s\, n}\, , e_{t\, n}$ associated with $\mu_n(s)\, , \mu_n(t)$, respectively, satisfy
\beq
\Bc(e_{s\, n}) \; \leq \; 0 \; \leq \Bc(e_{t\, n})
\eeq
\elem

\bpf
Let $\lambda_i\in (s, t)$ be a sequence monotonically increasing to $t$.
Using an argument similar to the proof of Lemma \ref{l.extract}, let $v_i\in V$ be a sequence of eigenvectors of $(a, b, m)$ associated with $\mu_n(\lambda_i)$ satisfying $v_i\rightarrow e_{t\, n}$ in $V$ as $i\rightarrow \infty$.
The orthogonality relation of Theorem \ref{t.ortho} then yields
\[
\lim_{i\rightarrow \infty}\frac{\mu_n(\lambda_i) \, - \, \mu_n(t)}{ \lambda_i \, - \, t} \; \; = \; \; \lim_{i\rightarrow \infty}\frac{-b(v_i, e_{t\, n})}{m(v_i, e_{t\, n})} \; \;  = \; \; -\, \frac{\Bc(e_{t\, n})}{\Mcl(e_{t\, n})}.
\]
Since $\mu_n(\lambda_i) >\mu_* = \mu_n(t)$ and $\lambda_i < t$, it follows that $\Bc(e_{t\, n})\geq 0$.  
The argument that $\Bc(e_{s\, n})\leq 0$ is similar.
\epf

If a component of $\mu_n(\lambda)$ is of the form $(-\infty, t)$ instead, then the proof above is exactly the same to show that an eigenvector $e_{t\, n}$ associated with $\mu_n(t)$ satisfies $\Bc(e_{t\, n}) \, \geq \, 0$ provided $\mu_*\notin \sigma(a, m)$.  
An analogous satement also holds for a component of $\mu_n(\lambda)$ of the form $(s, \infty)$.

Our main result for higher eigencurves is the following. 

\btm\label{t.gnbumps}
Assume $a, b, m$ satisfy $(A1)$-$(A3)$, and let $\mu_*\notin \sigma(a, m)$ be fixed.
Suppose $(s_1, t_1), \ldots, (s_n, t_n)$ are finite intervals with $(s_i, t_i)$ a component of some (not necessarily the same) $\{\mu_j>\mu_*\}$.
Two such intervals are either nested, with either endpoint possibly in common, or disjoint.  
If $(s_i, t_i) = (s_j, t_j)$, then these intervals are assumed to be components of two distinct eigencurves.
Then 
\beq
\mu_n(\lambda)\; \leq \; \mu_* \qquad \text{ for all }\quad \lambda \in (-\infty, \hat{s}\,] \cup [\, \hat{t}, \infty)
\eeq
where \; $\hat{s}:=\min\{s_1, \ldots, s_n\}$ \; and \; $\hat{t}:= \max\{t_1, \ldots, t_n\}$.
\etm

\bpf
For each $t_i$, with $1\leq i\leq n$, use Lemma \ref{l.eprobe} to get eigenvectors $e_{i\, n}$ associated with $\mu_n(t_i)$ satisfying $\Bc(e_{i\, n})\geq 0$, and define $Z:=\text{span}\{e_{1\, n}, \ldots, e_{n\,n}\}$.
By Theorem \ref{t.independent}, the vectors $e_{i\, n}$ are linearly independent so that $Z$ is an $n$-dimensional subspace of $V$.
By Corollary \ref{c.orthoH}, the vectors $e_{i\, n}$ are $b$-orthogonal as well as $(a-\mu_*m)$-orthogonal.
Let $v=\sum_{i=1}^n\alpha_i\,e_{i\,n}$, with $\alpha_i\in \mathbb{R}$, be a vector in $Z$ and let $\lambda \geq \hat{t}$.
Using the eigenpoint equation that each $e_{i\, n}$ satisfies and the orthogonality just mentioned yields the following:
\[
\Ac(v)  -  \mu_*\, \Mcl(v)  -  \lambda\, \Bc(v) 
\, =  \, \sum_{i=1}^n\alpha_i^2\, t_i\, \Bc(e_{i\, n}) \, - \, \lambda \sum_{i=1}^n\alpha_i^2\Bc(e_{i\, n}) 
\, =  \, \sum_{i=1}^n\alpha_i^2(t_i - \lambda)\, \Bc(e_{i\, n}) \; \leq  \; 0.
\]
Rearranging the exteme sides of this relation gives
\[
\mu_n(\lambda) \; \; \leq \; \; \sup_{v\in Z\setminus\{0\}} \left\{\frac{\Ac(v) - \lambda \Bc(v)}{\Mcl(v)} \right\} \; \; \leq \; \; \mu_* \qquad \text{for all }\lambda\geq \hat{t}
\]
using the variational characterization \eqref{e.minmax} for $\mu_n(\lambda)$.
A similar argument using  left endpoints $s_i$ shows that $\mu_n(\lambda)\leq \mu_*$ for all $\lambda \leq \hat{s}$ so that the desired result holds.
\epf

For fixed eigencurve $\mu_n(\lambda)$ and $\mu_*\in \R$ denote by $K_n(\mu_*)$ the number of components comprising the superlevel set $\{\mu_n > \mu_*\}$.
Lemma \ref{l.curve1} may then be expressed as saying that $K_1(\mu_*)\leq 1$ for all $\mu_*\in \R$.

The following result is an immediate consequence of Theorem \ref{t.gnbumps}.

\bcor\label{l.nbumps}
Assume $a, b, m$ satisfy $(A1)$-$(A3)$ and let $\mu_*\notin \sigma(a, m)$ be fixed. 
Then for each $n\in \mathbb{N}$ we have $K_n(\mu_*)\; \leq \; n$.
\ecor

\bpf 
Suppose $(s_1, t_1), \ldots, (s_n, t_n), (s_{n+1}, t_{n+1})$ are distinct components comprising the set $\{\mu_n > \mu_*\}$ with the endpoints satisfying $t_{i-1} < s_i < t_i < s_{i+1}$ for each $i$.
We point out that we could have $s_1 =-\infty$ and also $t_{n+1}=+\infty$.
Nevertheless, consider only the first $n$ components.
We conclude from the proof of Theorem \ref{t.gnbumps} that $\mu_n(\lambda)\leq \mu_*$ for $\lambda\geq t_n$ in contradiction to $(s_{n+1}, t_{n+1})$ being a component of $\mu_n(\lambda)$. 
Therefore, $\{\mu_n >\mu_*\}$ can have at most $n$ components so that $K_n(\mu_*)\leq n$ holds for each $n\in \mathbb{N}$.
\epf

Our main geometrical result for higher eigencurves is the following.

\btm
Assume $a, b, m$ satisfy $(A1)$-$(A3)$ and let $\mu_*\notin \sigma(a, m)$ be fixed.  Then for each $n\in\N$ we have
\beq
\sum_{i=1}^n K_i(\mu_*) \; \; \leq \;  \; n.
\eeq
\etm

\bpf
Let $(s_1, t_1), \ldots, (s_k, t_k)$ be distinct intervals, each a component of (not necessarily the same) $\{\mu_i > \mu_*\}$ for some $1\leq i\leq n$.
Define 
\[
\lambda_* \; := \; \inf\{s_i \; : \; 1\leq i\leq k\}\qquad \text{and}\qquad \lambda^* \; := \; \sup\{t_i \; : 1\leq i\leq k\}.
\]
Assume without loss of generality that $(s_1, t_1), \ldots, (s_{k_0}, t_{k_0})$ are all the components satisfying $t_i \leq t_{k_0} < \lambda^*$ for each $i=1, \ldots, k_0$. 
It follows from Theorem \ref{t.gnbumps} that 
\[
\mu_{k_0}(\lambda)\leq \mu_*\qquad \text{for all }\lambda\geq t_{k_0}.
\]
Since $\mu_i(\lambda)\leq \mu_{k_0}(\lambda)$ for each $i=1, \ldots, k_0$, it follows that $\mu_i(\lambda)\leq \mu_*$ for all $\lambda\geq t_{k_0}$ for each $i\leq k_0$.
For each $k_0< i\leq k$, we have that $t_i=\lambda^*$ and that $(s_i, t_i)$ is a component of one of the following
\[
\{\mu_{k_0+1}>\mu_*\}\, , \qquad \{\mu_{k_0+2}>\mu_*\}\, , \quad \ldots\quad , \{\mu_k >\mu_*\}.
\]
Two components $(s_i, t_i), (s_j, t_j)$ with $k_0 < i, j\leq k$ therefore belong to two different eigencurves above $\mu_*$.
Since the intervals considered here are components of $\{\mu_i>\mu_*\}$ for some $1\leq i\leq n$, it follows that $k- k_0 \leq n- k_0$, so that $k\leq n$.  
This means $\sum_{i=1}^nK_i(\mu_*) \; \leq \; n$ holds for each $n\in \mathbb{N}$.
\epf

The above analysis has dealt with the case of a horizontal line $\mu = \mu_*$ intersecting the spectrum $\sigma(a, b, m)$.
Consider now the case of a non-horizontal line $\mu = \alpha \, \lambda + \beta$ intersecting the spectrum of $(a, b, m)$.
The previous work generalizes to this case  using a simple change-of-variables as the relation
\[
a(e, v) \; = \;  \lambda \, b(e, v) \, + \mu\, m(e, v) \qquad \text{for all }v\in V,
\]
holds if and only if the following relation holds,
\[
a(e, v) \; = \; \lambda\, b_{\alpha}(e, v) \; + \; \beta\, m(e, v) \qquad \text{for all }v\in V
\]
where $b_{\alpha} : = b+ \alpha\, m$.
In this case, the $(a, b_{\alpha}, m)$-eigenproblem satisfies the same assumptions as the original $(a, b, m)$-eigenproblem.
Denote the eigencurves associated with the triple $(a, b_{\alpha}, m)$ by $\beta_{\alpha, n}(\lambda)$. 
The variational characterization \eqref{e.minmax} shows that
\[
\beta_{\alpha, n}(\lambda) \;  = \; \inf_{Z_k} \sup_{\substack{v\in Z_k \\ v\neq 0}} \left\{ \frac{\Ac(v) - \lambda\, \Bc_{\alpha}(v) }{\Mcl(v)}\right\} 
\;  =  \; \inf_{Z_k} \sup_{\substack{v\in Z_k \\ v\neq 0}} \left\{ \frac{\Ac(v) - \lambda\, \Bc(v) }{\Mcl(v)} \, - \, \alpha\, \lambda\right\} 
\;  =  \; \mu_n(\lambda) \, - \, \alpha \, \lambda
\]
where the minimization is over all $k$-dimensional subspaces $Z_k$ of $V$.
Thus, for fixed $\lambda\in \R$, the eigenvalues $\mu_n(\lambda)$ and $\beta_{\alpha, n}(\lambda)$ only differ by the constant $\alpha \, \lambda$, and therefore the eigencurves $(\lambda, \mu_n(\lambda))$ and $(\lambda, \beta_{\alpha, n}(\lambda))$ identify the same curve.
Hence, the results for horizontal-line-intersections generalizes.






\end{document}